\documentclass[titlepage,12pt]{article}
\usepackage{amssymb}
\usepackage{amsfonts}
\usepackage{amsmath}
\begin{document}

{\bf \Large Four Departures in Mathematics and Physics} \\ \\

{\bf Elemer E Rosinger} \\
Department of Mathematics \\
and Applied Mathematics \\
University of Pretoria \\
Pretoria \\
0002 South Africa \\
eerosinger@hotmail.com \\ \\

\hfill {\it Dedicated to Marie-Louise Nykamp} \\ \\

{\bf Abstract} \\

Much of Mathematics, and therefore Physics as well, have been limited by four rather consequential restrictions. Two of them are ancient taboos, one is an ancient and no longer felt as such bondage, and the fourth is a surprising omission in Algebra. The paper brings to the attention of those interested these four restrictions, as well as the fact that each of them has by now ways, even if hardly yet known ones, to overcome them. \\ \\

{\bf 1. Two Taboos and a Bondage form Ancient Times, \\
        \hspace*{0.5cm} plus a Modern Omission} \\

Logic, as is well known, is at the basis of Mathematics. Two {\it
ancient taboos} in Logic are :

\begin{itemize}

\item no contradictions

\item no self-reference

\end{itemize}

The reason for the first taboo seems so obvious as not to require
any argument. The reason for the second taboo is known at least
since ancient Greece and the Paradox of the Liar, with its modern
version being Russell's paradox in Set Theory. \\

In Mathematics, the {\it ancient bondage} is the tacit acceptance of
the Archimedean Axiom in Euclidean Geometry. The consequence is that the geometric line is coordinatized, ever since Descartes, by the
usual field $\mathbb{R}$ of real numbers. And as is well known, this is the only linearly
ordered
complete field that satisfies the Archimedean Axiom. \\
Amusingly, this uniqueness result has often been interpreted as an impressive justification of the identification between the geometric line and $\mathbb{R}$, thus preempting any considerations for other possibilities of coordinatization of the geometric line. \\
Consequently, the elementary fact is completely missed that this uniqueness does not happen because of certain seemingly fortunate mysterious reasons, but it is merely a direct consequence of the Archimedean Axiom. \\

As for the {\it modern omission}, it is in disregarding the fact that the
n-dimensional Euclidean space $\mathbb{R}^n$ is not only a module
over the field $\mathbb{R}$ of real numbers, but it has a
considerably {\it more rich structure} of module over the
non-commutative algebra $\mathbb{M}^n$ of n $\times$ n matrices of
real numbers. \\ \\

{\bf 2. No More ''No-Contradictions'' ...} \\

As seen in [1] and the literature cited there, recently, logical systems which are {\it inconsistent} have been introduced and studied. In this regard, quite unobserved by anybody, we have during the last more than six decades been massively using and relying upon such systems, ever since the introduction of digital electronic computers. Indeed, such computers, when seen as dealing with integer numbers, are based on the Peano Aximoms, plus the following one :

\begin{itemize}

\item There exists a positive integer $M >> 1$, such that $M + 1 = M$.

\end{itemize}

This $M$, usually larger than $10^{100}$, is called ''machine infinity'' and, depending on the particular computer, is the largest positive integer which that computer can rigorously represent. \\

Clearly, however, the above axiom is inconsistent with the Peano Axioms. \\

Nevertheless, as the more than six decades long considerable experience in using digital computers shows it, we can quite safely use and rely upon the trivially inconsistent axiom system upon which all digital computers are inevitably built. \\ \\

{\bf 3. No More ''No-Self-Reference'' ...} \\

The issue of self-reference is far more subtle and involved, and its roots appear to go back to prehistoric times. Indeed, anthropologists argue that three of the fundamental issues preoccupying humanity from time immemorial are change, infinity and self-reference. As for the latter, its importance may further be illustrated by the fact that in the Old Testament, in Exodus 3:14, it stands for no less than the self-proclaimed name of God which is ''I Am that I Am'' ... \\

A brief introduction, as well as the most important related literature can be found in [6], while an instructive simple application is presented in [5]. \\

{\bf 4. Non-Archimedean Mathematics : Reduced Power \\
        \hspace*{0.5cm} Algebras and Fields} \\

The way space, and also time, are - and have ever since ancient Egyptian Geometry as formalized by Euclid been - modelled mathematically is essentially centered on the usual common day to day human perspective. This, simplified to the one dimension of the geometric line, means among others that from any given point on the line one can reach every other point on that line after a finite number of steps of an a priori given size. It also means that, in one of the directions along such a line, there is one and only one ''infinity'' which in some rather strange way is not quite on the line, but as if beyond it. And a similar, second ''infinity'' lurks in the opposite direction. Furthermore, there cannot be smaller and smaller intervals of infinitesimal, yet nonzero length, just as much as there cannot be intervals with larger and larger infinite length. \\

On the other hand, modern Physics, with the advent of Cosmology based on General Relativity, and with the development of Quantum Field Theory, keeps massively - even if less obviously, and rather subtly - challenging the usual human perspective of space-time. And one of the consequences of not dealing seriously and systematically enough with that challenge is the presence of all sorts of so called ''infinities in Physics'', among others in Quantum Field Theory. Also, the human perspective based structure of space-time cannot find ''space'', or for that matter ''time'', where the so called multiverses may happen to exist, in case the Everett, Tegmark, and similar views of the quantum world or cosmos at large would indeed prove to have sufficient merit. \\

As for the possible merit of the idea of multiverse in general, let us recall the following estimates, [15, p.76]. It is considered that the visible part of Cosmos contains less than $10^{80}$ atoms. On the other hand, practically useful quantum algorithms, such as for instance, for factorization by Shor's method, may require quantum registers with at least $10^{500}$ states. And each such state must exist somewhere, since it can be involved in superpositions, entanglements, etc., which occur during the implementation of the respective quantum algorithm. \\

And then the question arises : \\

When such a quantum algorithm is effectively performed in a world with not more than $10^{80}$ atoms, then where does such a performance involving more than $10^{500}$ states happen ? \\

And yet amusingly, the ancient Egyptian thinking about the structure of space-time is so long and so hard entrenched that no one ever seems to wonder whether, indeed, there may not be other space-time mathematical structures that would be more appropriate for modern Physics. \\

Fortunately, as it happens, from mathematical point of view the situation can simply and clearly be described : the usual common day to day human perspective on the structure of space-time comes down mostly to nothing else but the acceptance of the Archimedean Axiom in Euclidean Geometry. \\

The surprisingly rich possibilities which open up once the Archimedean Axiom is set aside have been presented in certain details in [2-4,7,8], and the literature cited there. \\ \\

{\bf 5. Far Richer Module Structures for Space-Time} \\

It is quite surprising that for so long the following simple yet highly consequential algebraic fact has mostly been overlooked although it can offer considerable theoretical advantages in several branches of Linear Mathematics, as well as in applications. Indeed, so far, it has been Linear System Theory the only notable exception where a full use of that algebraic fact could be noted, [12]. \\

Let us start by considering any $n \geq 2$ dimensional vector space $E$ over some field $\mathbb{K}$. This means that the abelian group $( E, + )$ is a {\it module} over the unital commutative ring structure of $\mathbb{K}$, and for convenience, we shall consider it as customary, namely, a left module \\

(5.1)~~~ $ \mathbb{K} \times E \ni ( a, x ) ~\longmapsto~ a x \in E $ \\

Let us now denote by \\

(5.2)~~~ $ \mathbb{M}\,^n_\mathbb{K} $ \\

the set of all $n \times n$ matrices with elements in $\mathbb{K}$, which is a noncommutative unital algebra. Given now a basis $e_1, \ldots , e_n$ in $E$, we can define the {\it module structure} \\

(5.3)~~~ $ \mathbb{M}\,^n_\mathbb{K} \times E \ni ( A, x )~\longmapsto~ A x \in E $ \\

by the usual multiplication of a matrix with a vector, namely, if $A = ( a_{i, j} )_{1 \leq i, j \leq n} \in \mathbb{M}\,^n_\mathbb{K}$ and $x = x_1 e_1 + \ldots + x_n e_n \in E$, with $x_1, \ldots , x_n \in \mathbb{K}$, then \\

(5.4)~~~ $ A x =
               ( a_{1,1} x_1 + \ldots + a_{1,n} x_n ) e_1 + \ldots + ( a_{n,1} x_1 + \ldots + a_{n,n} x_n ) e_n  $ \\

Let us at last recall the existence of the {\it algebra embedding} \\

(5.5)~~~ $ \mathbb{K} \ni a ~\longmapsto~ I_a \in \mathbb{M}\,^n_\mathbb{K} $ \\

where $I_a$ is the $n \times n$ diagonal matrix with $a$ as the diagonal elements. \\

Obviously, for the one dimensional case $n = 1$, the module structures (5.1) and (5.3) are identical, and therefore, we shall not consider that case. \\

Two facts are important to note here, namely

\begin{itemize}

\item In view of the algebra embedding (5.5), the module structure (5.3) is obviously an {\it extension} of the original module structure (5.1).

\item Since for larger $n$, the algebra $\mathbb{M}\,^n_\mathbb{K}$ is considerably larger than the field $\mathbb{K}$, it follows that the module structure (5.3) is considerably more rich than the vector space structure (5.1).

\end{itemize}

The above obviously can be reformulated in a coordinate free manner for arbitrary finite or infinite dimensional vector spaces $E$ over a given field $\mathbb{K}$. Indeed, let us denote by ${\cal L}_\mathbb{K} ( E )$ the unital noncommutative algebra of $\mathbb{K}$-linear operators on $E$. Then we have the {\it algebra embedding} \\

(5.6)~~~ $ \mathbb{K} \ni a ~\longmapsto~ I_a \in {\cal L}_\mathbb{K} ( E ) $ \\

where $I_a ( x ) = a x$, for $x \in E$. \\

On the other hand, by the definition of ${\cal L}_\mathbb{K} ( E )$, we have the {\it module structure} \\

(5.7)~~~ $ {\cal L}_\mathbb{K} ( E ) \times E \ni ( A, x ) ~\longmapsto~ A x \in E $ \\

And then again, in view of (5.6), the module structure (5.7) is considerably more rich than the vector space structure
(5.1). \\

{\bf Remark 5.1.} \\

In certain ways, the module structures (5.1) and (5.7) on a given vector space $E$ are {\it extreme}. Indeed, in the
latter, $E$ can in fact be a {\it cyclic} module, [9-14], being generated by one single element. Namely, it is obvious
that, in case for instance $\mathbb{K}$ is the field of real or complex numbers, then the following holds \\

$~~~~~~ \exists~~ x \in E ~:~ \forall~~ y \in E ~:~ \exists~~ A \in {\cal L}_\mathbb{K} ( E ) ~:~ y = A x $ \\

Therefore, the issue is to find module structures on $E$ which in suitably ways are {\it intermediary} module structures
between those in (5.1) and (5.7).

\hfill $\Box$ \\

And therefore, one can ask :

\begin{itemize}

\item Which may possibly be the intermediary module structures between those in (5.1) and (5.7), and what can usefully be
done with such module structures on finite or infinite dimensional vector spaces ?

\end{itemize}

In this regard, let us first briefly recall what is known related to module structures defined by rings, and not merely by fields. Further details in this regard can be found in [9-11], while for convenience, a brief schematic account is presented in Appendix 1. \\
Here it should be mentioned that by far most of what is known in this regard concerns the finite dimensional cases, more precisely, the finitely generated modules over principal ideal domains, such as for instance presented in [9-11]. \\

However, even these simpler results have for more than four decades by now had fundamentally important applications in Linear System Theory, as originated and developed by R E Kalman, see [12, chap.10] and the references cited there. \\

Amusingly in this regard, back in 1969, in [12], Kalman found it necessary to object to the fact that theory of modules over principal ideal domains had not yet penetrated the university curriculum, and attributed that failure to possible pedagogical constraints. As it happens, however, the same situation keeps occurring more than four decades later, and the reasons for it are not any more clear now than were decades ago ... \\
By the way, Kalman also mentions in [12] that as far back as in 1931, Van der Waerden, in [13], which was to become the foundational textbook for Modern Algebra, had consistently used module theory in developing Linear Algebra, thus obtaining an impressive presentation of canonical forms and other important related results. \\
As for Bourbaki, [14], such a module theory based study of Linear Algebra was ever after to endorse in a significant measure, and also expand, the pioneering work in [13]. \\

Nowadays, the fast and massively developing theory of Quantum Computation has as one of its crucial and not yet solved problems the measurement of the extent of entanglement, entanglement being widely considered as perhaps the single most important resource upon which the performance of quantum computers can a does so much outperform that of usual electronic digital ones. And the mathematical background within which quantum entanglement occurs is the tensor product of complex Hilbert spaces which contain the relevant qubits involved in entanglement. \\

So far, however, when studying quantum entanglement, the only structure considered on the Hilbert spaces $H$ involved in the respective tensor products has been that of vector spaces over the field $\mathbb{C}$ of complex numbers, and not also the far richer structure of modules over the unital noncommutative algebras ${\cal L} ( H  )$ of linear bounded operators on $H$. \\

As for how more rich that latter structure is, we can recall that for a practically relevant quantum computation one is supposed to deal with quantum registers having at least $n > 100$ qubits. Consequently, the respective Hilbert spaces $H$ resulting from tensor products may have dimensions larger than $2^{100}$, with the corresponding yet far larger dimensions of the algebras ${\cal L} ( H  )$. \\

Certain indications about the impact of studying tensor products, and thus entanglement, based on the mentioned richer structure of Hilbert spaces are presented in section 6 in the sequel. \\

Needless to say, the issue of studying tensor products of vector spaces with the respective vector spaces considered endowed with the far richer module structures mentioned above can as well have an interest beyond Quantum Computation, and also within the general theory of Tensor Products, a theory which, being well known for its considerable difficulties, may possibly benefit from the introduction of such richer module structures, structures which in fact have always been there. \\
The number of major an long standing open problems in Linear Mathematics is well known, one of them being can the problem of invariant subspaces in Hilbert spaces. And again, the solution of this problem may possibly benefit from the consideration of Hilbert spaces $H$ not as mere vector spaces over the field $\mathbb{C}$ of complex numbers, but also as modules over the unital noncommutative algebra ${\cal L} ( H  )$ of linear bounded operators on $H$. \\

Let us now briefly review some of the relevant properties of modules over principal ideal domains. We recall that, for instance in  [10, p. 10], right from the very beginning the importance of the structure theorem of finitely generated modules over principal ideal domains is emphasized due, among others, to its far reaching consequences. \\

Let $R$ be a unital commutative ring and $M$ an $R$-module. Here and in the sequel, we shall mean left $R$-module whenever we write simply $R$-module. \\

A vector space $M$ is a particular cases of $R$-module, namely when $R$ is a field. And in case a vector space $M$ over the field $R$ has a finite dimension $1 \leq n < \infty$, then it has the following crucial property regarding its structure \\

(5.8)~~~ $ M \backsimeq R^n = R \bigoplus \ldots \bigoplus R $ \\

where the direct sum $\bigoplus$ has $n$ terms. \\

The point in the structure result (5.8) is the following. If $R$ is a field, then $R$ itself is the smallest nontrivial
$R$-module. Thus (5.8) means that every $n$-dimensional $R$-module $M$ is in fact the Cartesian product, or equivalently,
the direct sum of precisely $n$ smallest nontrivial $R$-modules. \\

Let us now see what happens when $R$ is not a field, but more generally, a unital commutative ring. A classical result which extends (5.8) holds when the following {\it two conditions} apply : \\

(5.9)~~~ $R$ is a principal ideal domain, or in short, PID \\

(5.10)~~~ $M$ is finitely generated over $R$ \\

The general structure result which extends (5.8) under the conditions (5.9), (5.10) is as follows, see Appendix 1 \\

(5.11)~~~ $ M \backsimeq R^s \bigoplus \, ( R / ( R a_1 ) ) \bigoplus \ldots \bigoplus \, ( R / ( R a_r ) ) $ \\

where $ a_1, \ldots , a_r \in R \setminus \{ 0, 1 \}$, with $a_i \, | \, a_{i + 1}$, for $1 \leq i \leq r - 1$, and the
above decomposition (5.11) is unique, up to isomorphism. \\

Again, the remarkable fact in the structure result (5.11) is the presence of the Cartesian product, or equivalently,
the direct sum $R^s$ of $s$ copies of $R$ itself, plus the presence of $r$ quotient $R$-modules $R / ( R a_i )$, where
each of the $R$-modules $R a_i $ is generated by one single element $a_i \in R$. \\

The case of interest for us will be the following particularization of the general module structures (5.7). \\

Given any field $\mathbb{K}$, let us consider the unital commutative algebra \\

(5.12)~~~ $ \mathbb{K} [ \xi ] $ \\

of polynomials in $\xi$ with coefficients in $\mathbb{K}$. As is well known, [9-14], this algebra is a principal ideal
domain, that is, PID. Furthermore, we have the algebra embedding \\

(5.13)~~~ $ \mathbb{K} \ni c \longmapsto c = c \xi^0 \in \mathbb{K} [ \xi ] $ \\

Now, for any given $A \in  {\cal L}_\mathbb{K} ( E )$, we consider the module structure on $E$, given by \\

(5.14)~~~ $ \mathbb{K} [ \xi ] \times E \ni ( \pi ( \xi ), x ) \longmapsto ( \pi ( A ) ) ( x ) \in E $ \\

that is, if $\pi ( \xi ) = \sum_{0 \leq i \leq n} \alpha_i \xi^i$, where $\alpha_i \in \mathbb{K}$, then $\pi ( A ) = \sum_{0 \leq i \leq n} \alpha_i A^i \in {\cal L}_\mathbb{K} ( E )$, thus indeed $( \pi ( A ) ) ( x ) \in E$. \\

It follows that, in (5.14), we have taken \\

(5.15)~~~ $ R = \mathbb{K} [ \xi ] $ \\

and \\

(5.16)~~~ $ M = E $ \\

hence (5.9), (5.10) are satisfied whenever $E$ is a finite dimensional vector space over $\mathbb{K}$. \\

With this particularization, instead of the general situation in (5.7), we shall obtain the family of module structures depending on $A \in  {\cal L}_\mathbb{K} ( E )$, namely \\

(5.17)~~~ $ \mathbb{K} [ \xi ] \times E \longrightarrow E $ \\

which obviously are still considerably more rich than the vector space structure (5.1), since the algebra $\mathbb{K}
[ \xi ]$ is so much more large than the field $\mathbb{K}$, see (5.13). \\

It is indeed important to note that each of the module structures (5.14), (5.17) depends on the specific linear operator
$A \in  {\cal L}_\mathbb{K} ( E )$, more precisely, these module structures are {\it not} the same for all such linear
operators. Therefore, instead of the given single vector space structure (5.1), we obtain the class of considerably richer
module structures (5.17). \\

And as seen in the next section, such a class of richer module structures can have an interest in the study of Quantum Computation,
related to a better understanding of tensor products of complex Hilbert spaces, and thus of the all important resource
given by entanglement. \\

We can conclude with the following clarification regarding the family of module structures (5.14), (5.17), [12, pp. 251,252 ] \\

{\bf Theorem 5.1.} \\

Let $E$ be a vector space over the field $\mathbb{K}$, and let $A \in  {\cal L}_\mathbb{K} ( E )$. Then $E$ admits the specific $\mathbb{K} [ \xi ]$ module structure (5.14), (5.17). \\

Conversely, any $\mathbb{K} [ \xi ]$ module structure on $E$ is of the form (5.14), (5.17), where $A \in  {\cal L}_\mathbb{K} ( E )$ is defined by \\

(5.18)~~~ $ E \ni x \longmapsto ( \pi_1 ( \xi ) ) ( x ) \in E $ \\

with the polynomial $\pi_1 ( \xi ) \in \mathbb{K} [ \xi ]$ given by $\pi_1 ( \xi ) = \xi$. \\

{\bf Remark 5.2.} \\

1) In Theorem 5.1. above the vector spaces $E$ can be arbitrary, thus need not be finite dimensional over $\mathbb{K}$, or
finitely generated over $\mathbb{K} [ \xi ]$. \\

2) Simple examples can show that decompositions (5.11) need not necessarily hold when the conditions (5.9), (5.10) are not
satisfied. \\
On the other hand, there is an obvious interest in more general decompositions of the form \\

(5.19)~~~ $ M \backsimeq R^I \bigoplus ( \bigoplus_{j \in J} \, ( R / ( R a_j ) ) ) $ \\

where $I$ or $J$, or both, may be infinite sets. \\
And clearly, the conditions (5.9), (5.10) are not necessary for such more general decompositions (5.19). Consequently,
there can be an interest in finding sufficient conditions for decompositions (5.19), conditions weaker than those in (5.9),
(5.10). \\ \\

{\bf 6. Entanglement and Modules over Principal Ideal Domains} \\

In view of its considerable interest in Quantum Computation, let us see what happens with the entanglement property of wave functions of composite quantum systems when the respective finite dimensional component complex Hilbert spaces are considered not merely as vector spaces over the field $\mathbb{C}$ of complex numbers, but also with the far richer structure of modules over principal ideal domains. \\

Specifically, let $H$ be an $n$-dimensional complex Hilbert space which is the state space of a given number of qubits describing a certain quantum computer register. For starters, we consider the larger quantum computer register with the double number of qubits, thus given by the usual tensor product $H \bigotimes H$. \\

As is known, [14], this tensor product is an $n^2$-dimensional Hilbert space which results as a particular case of the following general algebraic construction. \\

Let $E, F$ be two complex vector spaces. In order to define the tensor product $E \bigotimes F$, we start by denoting
with  \\

(6.1)~~~ $ Z $ \\

the set of all finite sequences of pairs \\

(6.2)~~~ $ ( x_1, y_1 ) \ldots ( x_m, y_m ) $ \\

where $m \geq 1$, while $x_1, \ldots , x_m \in E, y_1, \ldots , y_m \in F$. We consider on $Z$ the usual concatenation
operation of the respective finite sequences of pairs, and for algebraic convenience, we denote it by $\gamma$. It follows
that, in fact, (6.2) can be written as \\

(6.3)~~~ $ ( x_1, y_1 ) \ldots ( x_m, y_m ) = ( x_1, y_1 ) \gamma \ldots \gamma ( x_m, y_m ) $ \\

that is, a concatenation of single pairs $( x_i, y_i ) \in E \times E$. Obviously, if at least one of the vector spaces
$E$ or $F$ is not trivial, thus it has more than one single element, then the above binary operation $\gamma : Z \times Z
\longrightarrow Z$ is {\it not} commutative. \\

Also, we have the injective mapping \\

(6.4)~~~ $ E \times F \ni ( x, y ) \longmapsto ( x, y ) \in Z $ \\

which will play an important role. \\

Now, one defines on $Z$ the {\it equivalence} relation $\approx$ as follows. Two finite sequences of pairs in $Z$ are equivalent, iff \\

(6.5)~~~ they are the same \\

or they can be obtained form one another by a finite number of applications of the following four operations \\

(6.6)~~~ permuting the pairs in the sequences \\

(6.7)~~~ replacing a pair $( x + x', y )$ with the pair $( x, y ) \gamma ( x', y )$ \\

(6.8)~~~ replacing a pair $( x, y + y' )$ with the pair $( x, y ) \gamma ( x, y' )$ \\

(6.9)~~~ replacing a pair $( c x, y )$ with the pair $( x , c y )$, where $c \in \mathbb{C}$ is \\
         \hspace*{1.3cm} any complex number \\

And now, the usual tensor product $E \bigotimes F$ is obtained simply as the {\it quotient} space \\

(6.10)~~~ $ E \bigotimes F = Z / \approx $ \\

According to usual convention, the coset, or the equivalence class in $ E \bigotimes F$ which corresponds to the element
$( x_1, y_1 ) \gamma \ldots \gamma ( x_m, y_m ) \in Z$ is denoted by \\

(6.11)~~~ $ ( x_1 \bigotimes y_1 ) + \ldots + ( x_m \bigotimes y_m ) \in E \bigotimes F $ \\

and thus in view of (6.4), one has the following embedding \\

(6.12)~~~ $ E \times F \ni ( x, y ) \longmapsto E \bigotimes F $ \\

We note that in view of (6.9), the equivalence relation $\approx$ on $Z$ depends on the field $\mathbb{C}$ of complex numbers, therefore, we shall denote it by $\approx_\mathbb{C}$. Consequently, the usual tensor product $\bigotimes$ in (6.10) also depends on $\mathbb{C}$. In view of that, it will be convenient to denote this usual tensor product by \\

(6.13)~~~ $ E \bigotimes_\mathbb{C} F $ \\

since in the sequel other tensor products are going to be considered. \\

Indeed, let, see Appendix 2 \\

(6.14)~~~ $ A \in {\cal L}_\mathbb{C} ( E ),~~~ B \in {\cal L}_\mathbb{C} ( F ) $ \\

and let us consider the respective $\mathbb{C} [ \xi ]$-module structures (5.14), (5.17) on $E$ and $F$. Then the
construction in (6.1) - (6.10) above can be repeated, with the only change being that, instead of (6.9), we
require the condition \\

(6.15)~~~ replacing a pair $( \pi (  \xi ) x, y )$ with the pair $( x , \pi ( \xi ) y )$, where \\
          \hspace*{1.4cm} $\pi ( \xi )\in \mathbb{C} [ \xi ]$ is any polynomial in $\xi$, with complex coefficients \\

Here we note that the operation $\pi (  \xi ) x$ takes place in the $\mathbb{C} [ \xi ]$-module $E$, and thus in view of
(5.14), (5.17), it depends on $A$. Similarly, the operation $\pi (  \xi ) y$ takes place in the $\mathbb{C}
[ \xi ]$-module $F$, and hence it depends on $B$. \\
Therefore, the resulting equivalence relation we denote by $\approx_{A, B}$. Consequently, the corresponding tensor
product is denoted by \\

(6.16)~~~ $ E \bigotimes_{A, B} F = Z / \approx_{A, B} $ \\

Regarding the relationship between the two tensor products (6.13) and (6.16), it is easy to see that, for $s, t \in Z$, we have \\

(6.17)~~~ $ s \approx_\mathbb{C} t ~~\Longrightarrow~~ s \approx_{A, B} t $ \\

Indeed, given $x \in E, y \in F, c \in \mathbb{C}$, then in view of (6.9), we have \\

$~~~~~~ ( c x, y ) \approx_\mathbb{C} ( x , c y ) $ \\

However, (5.13) gives $c = c \xi^0 \in \mathbb{C} [ \xi ]$, thus (6.15) results in \\

$~~~~~~ ( c x, y ) \approx_{A, B} ( x , c y ) $ \\

And obviously, $( \pi ( \xi ) x, y ) \approx_{A, B} ( x , \pi ( \xi ) y )$, for some $\pi ( \xi ) \in \mathbb{C} [ \xi ]$,
need not, in view of (5.13), always imply any equivalence $\approx_\mathbb{C}$. \\

Consequently, for each $s \in Z$, the cosets, or the equivalence classes $( s )_{\approx_\mathbb{C}} \in E
\bigotimes_\mathbb{C} F$ and $( s )_{\approx_{A, B}} \in E \bigotimes_{A, B} F$ are in the relation \\

(6.18)~~~ $ ( s )_{\approx_\mathbb{C}} \subseteq ( s )_{\approx_{A, B}} $ \\

thus we have the following {\it surjective} linear mapping \\

(6.19)~~~ $ E \bigotimes_\mathbb{C} F \ni ( s )_{\approx_\mathbb{C}} ~~\longmapsto~~
                                ( s )_{\approx_{A, B}} \in E \bigotimes_{A, B} F $ \\

{\bf Remark 6.1.} \\

1) In the right hand side of (6.19) we have a large family of tensor products corresponding to various $A \in
{\cal L}_\mathbb{C} ( E )$, and $B \in {\cal L}_\mathbb{C} ( F )$. In this regard we can note that, even if $E = F$, one
can still choose $A \neq B$. \\

2) In particular, we recall once more that, in Quantum Computation, practically relevant quantum registers are supposed to
be modelled by the tensor product of $n > 100$ single qubit registers, each of them described by the complex Hilbert space
$H = \mathbb{C}^2$, and thus resulting in the $n$-factor usual tensor product \\

(6.20)~~~ $ H \bigotimes_\mathbb{C} \ldots \bigotimes_\mathbb{C} H = \mathbb{C}^{( 2^n )} $ \\

Also, we can again recall that {\it entanglement} is considered to be one of the most important resources in making
quantum computation so much more powerful than that with usual digital electronic computers. \\

And in order to try to grasp the magnitude of the problems related to entanglement, we can note that the version of (6.12)
which corresponds to (6.20) is the injective mapping \\

(6.21)~~~ $ H \times \ldots \times H \ni ( x_1, \ldots , x_n ) \longmapsto
               x_1 \bigotimes_\mathbb{C} \ldots \bigotimes_\mathbb{C} x_n \in
                          H \bigotimes_\mathbb{C} \ldots \bigotimes_\mathbb{C} H $ \\

where in the left hand side we have a $2 n$-dimensional complex Hilbert space, while the Hilbert space in right hand side
has $2^n$ dimensions, and here the cases of interest are when we are dealing with $n > 100$. \\

Consequently, the amount of {\it unentangled} elements in the tensor product in (6.20), namely those given by $x_1
\bigotimes_\mathbb{C} \ldots \bigotimes_\mathbb{C} x_n$ in (6.21), are immensely smaller than the rest of the elements in
the tensor product (6.20), elements which by definition are {\it entangled}. \\

And then, the interest in surjective mappings (6.19) comes from the fact that, with the large variety of the other tensor
products $\bigotimes_{A, B}$ which they involve in addition to the usual one $\bigotimes_\mathbb{C}$, one may obtain
further useful and not yet studied possibilities in order to better understand the phenomenon of entanglement. \\

3) Certainly, so far, Quantum Mechanics, as a theory of Physics, has only used the tensor product $\bigotimes_\mathbb{C}$
when modelling the composition of quantum systems. \\

However, even if the tensor products $\bigotimes_{A, B}$ may never be found to directly model effective physical quantum
processes, they may still turn out to have an advantage as a mathematical instrument in studying composite quantum
systems. \\
In this regard, and as mentioned at pct. 1) above, it may possibly be of interest to study the composite quantum system in
(6.20) by its various surjective mappings into tensor products of type $\bigotimes_{A, B}$, with each of the $n$ factors
$H$ having its own $A_i \in {\cal L}_\mathbb{C} ( H )$ used in the definition of the respective tensor products. \\

{\bf Example 6.1.} \\

Let us consider the simple situation of a quantum register with two qubits. Their state space is given by the usual tensor product \\

(6.22)~~~ $ H \bigotimes_\mathbb{C} H = H \bigotimes H = \mathbb{C}^2 \bigotimes \mathbb{C}^2 = \mathbb{C}^4 $ \\

Let us take two linear operators \\

(6.23)~~~ $ A, B \in {\cal L}_\mathbb{C} ( \mathbb{C}^2 ) $ \\

and consider the corresponding tensor products, see (6.16) \\

(6.24)~~~ $ \mathbb{C}^2 \bigotimes_{A, B} \mathbb{C}^2 $ \\

with the associated surjective linear mappings, see (6.19) \\

(6.25)~~~ $ \mathbb{C}^2 \bigotimes_\mathbb{C} \mathbb{C}^2 \ni ( s )_{\approx_\mathbb{C}} \longmapsto
                               ( s )_{\approx_{A, B}} \in  \mathbb{C}^2 \bigotimes_{A, B} \mathbb{C}^2 $ \\

We shall consider particular instances of $A$ and $B$ in (6.23), and establish the corresponding specific forms of the tensor products $\mathbb{C}^2 \bigotimes_{A, B} \mathbb{C}^2$ in (6.24) and (6.25). \\

Let us take specifically \\

(6.26)~~~ $ A = \left ( \begin{array}{l} a ~~ 0 \\
                                         0 ~~ b
                        \end{array} \right ),~~~~  B = \left ( \begin{array}{l} c ~~ 0 \\
                                                                                0 ~~ d
                                                               \end{array} \right ) $ \\

where $a, b, c, d \in \mathbb{C}$ are certain given complex numbers. \\

We establish now the $\mathbb{C} [ \xi ]$-module structure on $\mathbb{C}^2$ which corresponds to $A$, and similarly, the $\mathbb{C} [ \xi ]$-module structure on $\mathbb{C}^2$ which corresponds to $B$. \\

In view of (5.14), we have with respect to the first mentioned module is \\

(6.27)~~~ $ \mathbb{K} [ \xi ] \times \mathbb{C}^2 \ni ( \pi ( \xi ), x ) \longmapsto
                                                    ( \pi ( A ) ) ( x ) \in \mathbb{C}^2 $ \\

where for $\pi ( \xi ) = \sum_{0 \leq i \leq n} \alpha_i \xi^i$ and $x = \left ( \begin{array}{l} u \\
                                                                                                  v
                                                                                 \end{array} \right )$, we have \\ \\

(6.28)~~~ $ ( \pi ( A ) ) ( x ) =  \left ( \begin{array}{l} u \pi ( a ) \\
                                                            v \pi ( b )
                                           \end{array} \right ) = \left ( \begin{array}{l}
                                                                              u \sum_{0 \leq i \leq n} \alpha_i a^i \\
                                                                              v \sum_{0 \leq i \leq n} \alpha_i b^i
                                                                           \end{array} \right ) $ \\

Similarly, the second mentioned module is \\

(6.29)~~~ $ \mathbb{K} [ \xi ] \times \mathbb{C}^2 \ni ( \pi ( \xi ), y ) \longmapsto
                                                    ( \pi ( B ) ) ( y ) \in \mathbb{C}^2 $ \\

where for $\pi ( \xi ) = \sum_{0 \leq i \leq n} \alpha_i \xi^i$ and $y = \left ( \begin{array}{l} w \\
                                                                                                  z
                                                                                 \end{array} \right )$, we have \\ \\

(6.30)~~~ $ ( \pi ( B ) ) ( y ) =  \left ( \begin{array}{l} w \pi ( c ) \\
                                                            z \pi ( d )
                                           \end{array} \right ) = \left ( \begin{array}{l}
                                                                              w \sum_{0 \leq i \leq n} \alpha_i c^i \\
                                                                              z \sum_{0 \leq i \leq n} \alpha_i d^i
                                                                           \end{array} \right ) $ \\

The effect of the above module actions (6.27), (6.29) on the tensor product $\mathbb{C}^2 \bigotimes_{A, B} \mathbb{C}^2$ appears only in the application of (6.15). Namely, we have for instance \\

(6.31)~~~ $ \left ( \left ( \begin{array}{l} u \pi ( a ) \\
                                             v \pi ( b )
                            \end{array} \right ) , \left ( \begin{array}{l} w \\
                                                                            z
                                                           \end{array} \right ) \right ) \approx_{A, B} 
             \left ( \left (  \begin{array}{l} u \\
                                               v   
                              \end{array} \right ) ,  \left ( \begin{array}{l} w \pi ( c ) \\
                                                                               z \pi ( d )
                                                              \end{array} \right ) \right ) $ \\
                                                              
for $u, v, w, z \in \mathbb{C}$ and $\pi ( \xi ) \in \mathbb{C} [ \xi ]$, and thus \\

(6.32)~~~ $ \left ( \begin{array}{l} u \pi ( a ) \\
                                     v \pi ( b )
                    \end{array} \right ) \bigotimes_{A, B} \left ( \begin{array}{l} w \\
                                                                                    z
                                                                   \end{array} \right ) =
            \left (  \begin{array}{l} u \\
                                      v
                      \end{array} \right ) \bigotimes_{A, B} \left ( \begin{array}{l} w \pi ( c ) \\
                                                                                      z \pi ( d )
                                                                     \end{array} \right ) $ \\

Now, through the mapping (6.25), we shall obtain \\

(6.33)~~~ $ \left ( \begin{array}{l} u \pi ( a ) \\
                                     v \pi ( b )
                    \end{array} \right ) \bigotimes_\mathbb{C} \left ( \begin{array}{l} w \\
                                                                                        z
                                                                       \end{array} \right ) \longmapsto 
             \left ( \begin{array}{l} u \pi ( a ) \\
                                      v \pi ( b )
                    \end{array} \right ) \bigotimes_{A, B} \left ( \begin{array}{l} w \\
                                                                                    z
                                                                   \end{array} \right ) $ \\ \\  

(6.34)~~~ $ \left (  \begin{array}{l} u \\
                                      v
                      \end{array} \right ) \bigotimes_\mathbb{C} \left ( \begin{array}{l} w \pi ( c ) \\
                                                                                          z \pi ( d )
                                                                         \end{array} \right ) \longmapsto 
            \left (  \begin{array}{l} u \\
                                      v
                      \end{array} \right ) \bigotimes_{A, B} \left ( \begin{array}{l} w \pi ( c ) \\
                                                                                      z \pi ( d )
                                                                     \end{array} \right ) $ \\
                                                                     
consequently, in view of (6.32), we obtain the simplification in the usual tensor product given by \\

(6.35)~~~ $ \left ( \begin{array}{l} u \pi ( a ) \\
                                     v \pi ( b )
                    \end{array} \right ) \bigotimes_\mathbb{C} \left ( \begin{array}{l} w \\
                                                                                        z
                                                                       \end{array} \right ) \longmapsto s $ \\ \\
                                                                       
(6.36)~~~ $ \left (  \begin{array}{l} u \\
                                      v
                      \end{array} \right ) \bigotimes_\mathbb{C} \left ( \begin{array}{l} w \pi ( c ) \\
                                                                                          z \pi ( d )
                                                                         \end{array} \right ) \longmapsto s $ \\ \\
                                                                         
Such possible simplifications of the usual tensor product form a large family, since they depend on the operators $A, B$ in (6.23) and on the polynomials $\pi ( \xi ) \in \mathbb{C} [ \xi ]$. \\

And such simplifications can therefore contribute to a better understanding of entanglement, with a respective study to be published elsewhere. \\ 

Another subsequent study is to clarify the simplifications of type (5.11) in the structure of quantum registers (6.20) when they are considered with their module structure over the principal ideal domain $\mathbb{C} [ \xi ]$. \\ \\ 

{\bf Appendix 1 : Structure of Modules on Principal Ideal \\
                  \hspace*{2.8cm} Domains} \\

( Brief summary following [9] ) \\

Basic assumptions \\

1) $R$ unital ring \\

2) $M$ left $R$-module \\ \\

(1) Definition, p. 272 \\

If $R$ ring, $M$ left $R$-module, $X \subseteq M$, define \\

$ < X > = \{ \sum_{i \in I} a_i x_i ~|~ a_i \in R, x_i \in X, I $ finite $ \}$ \\

If $M$ finitely generated, define \\

$ rank_R ( M ) = \inf \{ car X ~|~ X \subseteq M, M = < X > \}  $ \\

Define $X$ linearly independent, iff for $x_1, \ldots , x_n \in X$, pairwise different, $a_1, \ldots , a_n \in R$, one has \\

$ \sum_{i \in I} a_i x_i = 0 ~~\Longrightarrow~~ a_1 = \ldots = a_n = 0 $ \\

Define $X$ basis for $M$, iff \\

$ M = < X >,~~ X $ linearly independent \\

Define $M$ free, iff $M$ has a basis $X$. Then call $M$ free on $X$. \\

(2) Proposition 1, p. 273 \\

If $M$ finitely generated free module over $R$, then all the bases of $M$ have the same cardinal. \\

(3) Theorem 3, p. 274 \\

Every $R$-module $M$ is the quotient of a free $R$-module. \\

(4) Definition, p. 276 \\

If $I$ is a left ideal in $R$, define $I$ principal, iff \\

$ \exists a \in R ~:~ I = R a $ \\

Define $R$ principal ideal ring, iff all left and right ideals are principal. \\

Define $R$ principal ideal domain, or PID, iff $R$ principal ideal ring and integral domain. \\

(5) Definition, p. 283 \\

If $R$ integral domain, $M$ $R$-module, $u \in M$, define $u$ torsion element, iff \\

$ \exists a \in R, a \neq 0 ~:~ a u = 0 $ \\

Define torsion submodule \\

$ t M = \{ u \in M ~|~ u $ torsion element $ \}$ \\

Define \\

$M$ torsion module $~~\Longleftrightarrow~~ M = t M $ \\

$M$ torsion free module $~~\Longleftrightarrow~~ t M = \{ 0 \} $ \\

(6) If $R$ integral domain, $u \in M$, then \\

a) $ R u $ free on $u ~~\Longleftrightarrow~~ u \notin t M $ \\

b) $ t M $ $R$-module \\

c) $ t t M = t M $, $ t M $ torsion module, $ M / t M $ torsion free module \\

(7) Proposition 1, p. 283 \\

If $R$ PID, $n \geq 1$, then every submodule $M \subseteq R^n$ is free, with $rank_R ( M ) \leq n$. \\

(8) Corollary 2, p. 284 \\

If $R$ PID, $M$ finitely generated $R$-module, $N \subseteq M$ submodule, then \\

a) $N$ finitely generated \\

b) $ rank_R ( N ) \leq rank_R ( M ) $ \\

{\bf (9) Theorem 2 : Decomposition}, p. 284 \\

If $R$ PID, $M$ finitely generated $R$-module, then \\

$ M \backsimeq R^s \bigoplus R / ( R a_1 ) \bigoplus \ldots \bigoplus R / ( R a_r ) $ \\

where \\

$ a_1, \ldots , a_r \in R \setminus \{ 0, 1 \} $ \\

$ a_i \, | \, a_{i + 1},~~~~ 1 \leq i \leq r - 1 $ \\

and the above is unique, up to isomorphism. \\

Also \\

$ s = rank_R ( M ) $ \\

$ a_1, \ldots , a_r \in R \setminus \{ 0, 1 \} $ invariant factors of $M$. \\

{\bf (10) Corollary}, p. 285 \\

If $R$ PID, $M$ finitely generated $R$-module, then \\

$ M $ free $~~\Longleftrightarrow~~ M$ torsion free \\

(11) Defintion, p. 285 \\

If $R$ PID, $p \in R$ prime, $M$ $R$-module, define \\

$ M $ $p$-primary module $~~\Longleftrightarrow~~
            M = M_p = \{ x \in M ~|~ \exists n \geq 0 ~:~ p^n x = 0 \} $ \\

$ M $ primary module $~~\Longleftrightarrow~~ \exists p \in R $ prime ~:~ $ M $ $p$-primary module \\

{\bf (12) Theorem 3 Decomposition}, p. 285 \\

If $R$ PID, $M$ torsion $R$-module, then \\

$ M = \bigoplus_{i \in I} M_{p_i} $ \\

with $p_i \in R$ pairwise different primes. \\

If $M$ is finitely generated, then $I$ finite. \\

The decomposition is unique. \\

{\bf (13) Corollary}, p. 286 \\

If $R$ PID, $M$ finitely generated $R$-module, then \\

$ M = \bigoplus R / ( R p_i^{\alpha_{i, j}}) $ \\

with $p_i \in R$ pairwise different primes, and \\

$ \alpha_{i, 1} \leq \ldots \leq \alpha_{i, r_i} $ \\ \\

{\bf Appendix 2 : Tensor Products} \\

{\bf The Usual Tensor Product.} Let $R$ be a unital ring and $E, F$ two left $R$-modules. Then the tensor product \\

(A2.1)~~~ $ E \bigotimes_R F $ \\

is defined as follows. We denote with  \\

(A2.2)~~~ $ Z $ \\

the set of all finite sequences of pairs \\

(A2.3)~~~ $ ( x_1, y_1 ) \ldots ( x_m, y_m ) $ \\

where $m \geq 1$, while $x_1, \ldots , x_m \in E, y_1, \ldots , y_m \in F$. Further, we consider on $Z$ the usual concatenation operation of the respective finite sequences of pairs, and for algebraic convenience, we denote it by $\gamma$. It follows that, in fact, (A2.3) can be written as \\

(A2.4)~~~ $ ( x_1, y_1 ) \ldots ( x_m, y_m ) = ( x_1, y_1 ) \gamma \ldots \gamma ( x_m, y_m ) $ \\

that is, a concatenation of single pairs $( x_i, y_i ) \in E \times F$. Obviously, if at least one of the modules $E$ or $F$ is not trivial, thus it has more than one single element, then the above binary operation $\gamma : Z \times Z
\longrightarrow Z$ is {\it not} commutative. \\

Also, we have the injective mapping \\

(A2.5)~~~ $ E \times F \ni ( x, y ) \longmapsto ( x, y ) \in Z $ \\

which will play an important role. \\

Now, we define on $Z$ the {\it equivalence} relation $\approx_R$ as follows. Two finite sequences of pairs in $Z$ are equivalent, iff \\

(A2.6)~~~ they are the same \\

or they can be obtained form one another by a finite number of applications of the following four operations \\

(A2.7)~~~ permuting the pairs in the sequences \\

(A2.8)~~~ replacing a pair $( x + x', y )$ with the pair $( x, y ) \gamma ( x', y )$ \\

(A2.9)~~~ replacing a pair $( x, y + y' )$ with the pair $( x, y ) \gamma ( x, y' )$ \\

(A2.10)~~~ replacing a pair $( c x, y )$ with the pair $( x , c y )$, where $c \in R$ \\

And now, the tensor product $E \bigotimes_R F$ is obtained simply as the {\it quotient} space \\

(A2.11)~~~ $ E \bigotimes_R F = Z / \approx_R $ \\

According to usual convention, the coset, or the equivalence class in $ E \bigotimes_R F$ which corresponds to the element $( x_1, y_1 ) \gamma \ldots \gamma ( x_m, y_m ) \in Z$ is denoted by \\

(A2.12)~~~ $ ( x_1 \bigotimes_R y_1 ) + \ldots + ( x_m \bigotimes_R y_m ) \in E \bigotimes_R F $ \\

and thus in view of (A2.5), one has the following embedding \\

(A2.13)~~~ $ E \times F \ni ( x, y ) \longmapsto E \bigotimes_R F $ \\

We note that in view of (A2.10), the equivalence relation $\approx_R$ on $Z$ may indeed depend on $R$. Consequently, the tensor product $\bigotimes_R$ in (A2.11) may also depend on $R$. \\

{\bf Tensor Product over a Subring.} Let us now consider $S \subseteq R$ a unital subring in $R$. Then clearly $E, F$ are also $S$-modules, and one can construct as above the corresponding tensor product \\

(A2.14)~~~ $ E \bigotimes_S F $ \\

Therefore, the question is : what is the relationship between the tensor products (A2.11) and (A2.14) ? \\

The answer is easy to establish by noting that in the definition of the tensor product (A2.14) by the corresponding equivalence relationship $\approx_S$, the above relations (A2.2) - (A2.9) remain the same, while (A2.10) is changed
to \\

(A2.15)~~~ replacing a pair $( c x, y )$ with the pair $( x , c y )$, where $c \in S$ \\

Consequently, for $s, t \in Z$, we have \\

(A2.16)~~~ $ s \approx_S t ~~\Longrightarrow~~ s \approx_R t $ \\

Indeed, if $( x, y ) \in Z, c \in S$, then (A2.15) gives $( c x, y ) \approx_S ( x, c y )$. On the other hand, $S \subseteq R$ and (A2.10) lead to $( c x, y ) \approx_R ( x, c y )$, since $c \in R$. However, we need not always have the converse, namely, when $c \in R \setminus S$, and thus (A2.10) yields $( c x, y ) \approx_R ( x, c y )$, it need not mean that (A2.15) applies, hence giving $( c x, y ) \approx_S ( x, c y )$, since we have now $c \notin S$. \\

And now in view of (A2.16), we obtain the natural surjective linear mapping \\

(A2.17)~~~ $ E \bigotimes_S F \ni ( s )_{\approx_S} \longmapsto ( s )_{\approx_R} \in E \bigotimes_R F $ \\

{\bf Two Branching Tensor Products.} Let $f : R \longrightarrow S, g : R \longrightarrow T$ be two ring homomorphisms. Further, let $E$ be an $S$-module, while $F$ is a $T$-module. \\

In defining an equivalence relation on $Z$ we keep again (A2.2) - (A2.9), while we replace (A2.10) with \\

(A2.18)~~~ replacing a pair $( f ( c ) x, y )$ with the pair $( x , g ( c ) y )$, where \\
           \hspace*{1.8cm} $c \in R$ \\

Clearly, the resulting equivalence relation on $Z$ may depend on $f$ and $g$, thus we denote it by
$\approx_{f, \, g}$. \\

And now, we can define the corresponding tensor product \\

(A2.19)~~~ $ E \bigotimes_{f, \, g} F = Z / \approx_{f, \, g} $ \\

{\bf Example} \\

Let $R = S = T = E = F = \mathbb{R}$, while $f = id_\mathbb{R}$, and $g ( c ) = a c$, for $c \in \mathbb{R}$, where $a \in \mathbb{R}$ is given. \\

If now $x, y \in \mathbb{R}$, and $c \in \mathbb{R}$, then (A2.18) implies that \\

(A2.20)~~~ $ ( c x, y ) \approx_{f, \, g} ( x, a c y ) $ \\


\begin{thebibliography} {99}

\bibitem{} Rosinger E E : On the Safe Use of Inconsistent Mathematics. arXiv:0811.2405

\bibitem{} Rosinger E E : Heisenberg Uncertainty in Reduced Power Algebras. arXiv:0901.4825

\bibitem{} Rosinger E E : No-Cloning in Reduced Power Algebras Authors. arXiv:0902.0264

\bibitem{} Rosinger E E : Special Relativity in Reduced Power Algebras. arXiv:0903.0296

\bibitem{} Rosinger E E, Van Zyl A : Self-Referential Definition of Orthogonality. arXiv:0904.0082

\bibitem{} Rosinger E E : Brief Lecture Notes on Self-Referential Mathematics, and Beyond. arXiv:0905.0227

\bibitem{} Rosinger E E : Where Infinitesimals Come From ... \\ arXiv:0909.4396

\bibitem{} Rosinger E E : Surprising Properties of Non-Archimedean Field Extensions of the Real Numbers. arXiv:0911.4824

\bibitem{} Cohn P M : Algebra, Vol. 1. Wiley, London, 1974

\bibitem{} Blyth T S : Module Theory, An Approach to Linear Algebra. Oxford Science Publications. Clarendon Press, Oxford, 1990

\bibitem{} Adkins W A, Weintraub S H : Algebra, An Approach via Module Theory. Springer, New York, 1992

\bibitem{} Kalman R E, Falb P L, Arbib M A : Topics in Mathematical System Theory. McGraw-Hill, New York, 1969

\bibitem{} Van der Waerden B L : Moderne Algebra. Springer, Heidelberg, 1931

\bibitem{} Bourbaki N : El\'{e}ments de Math\'{e}matique, Livre II, Alg\'{e}bre. Hermann, Paris, 1962

\bibitem{} Gribbin, J : In Search of the Multiverse. Penguin Books, London, 2009

\end{thebibliography}
\end{document}